# Deep Learning-based Model Predictive Control for Resonant Power Converters

S. Lucia, *Member, IEEE*, D. Navarro, B. Karg, H. Sarnago, *Senior Member, IEEE*, and O. Lucía, *Senior Member, IEEE*.

*Abstract-* Resonant power converters offer improved levels of efficiency and power density. In order to implement such systems, advanced control techniques are required to take the most of the power converter. In this context, model predictive control arises as a powerful tool that is able to consider nonlinearities and constraints, but it requires the solution of complex optimization problems or strong simplifying assumptions that hinder its application in real situations. Motivated by recent theoretical advances in the field of deep learning, this paper proposes to learn, offline, the optimal control policy defined by a complex model predictive formulation using deep neural networks so that the online use of the learned controller requires only the evaluation of a neural network. The obtained learned controller can be executed very rapidly on embedded hardware. We show the potential of the presented approach on a Hardware-in-the-Loop setup of an FPGA-controlled resonant power converter.

## I. INTRODUCTION

Since the appearing of resonant power conversion in the 70s [1-3], this technology has significantly evolved to provide high performance and cost-effective power conversion in many different applications. Nowadays, the combination of this technology with advances in digital control and semiconductor technology [4] has enabled the design of power electronic systems with a superior performance grade.

The design of such systems also requires advanced control techniques to take the most of these power converters [5]. This includes also the development of modern digital control architectures [6], based either on FPGA [7, 8] or microprocessor/DSP, that are enhanced by new development techniques such as high level synthesis [9-11] and automatic code generation [12], enabling the design of more advanced and optimized controller implementations.

In recent years, a number of new control techniques have arisen including sliding-mode control [13] and model-free control [14, 15] applied to MIMO systems. In this context, model predictive control (MPC) has also arisen as an effective tool to control complex systems [16] even when fast sampling times are required. MPC uses a mathematical model that describes the system dynamics to predict the future trajectories of the system states and to compute a sequence of

Manuscript received July 23, 2019; revised October 14, 2019; accepted January 22, 2020. Copyright © 2020 IEEE. Personal use of this material is permitted. However, permission to use this material for any other purposes must be obtained from the IEEE by sending a request to pubs-permissions@ieee.org.
This work was partly supported by the Spanish MINECO under Project TEC2016-78358-R, by the Spanish MICINN and AEI under Project RTC-2017-5965-6, co-funded by EU through FEDER program, by the DGA-FSE, and by DGA under Project LMP106-18.
O. Lucia, D. Navarro and H. Sarnago are with the Department of Electronic Engineering and Communications, University of Zaragoza, SPAIN (phone: +34976761000, e-mail: olucia@unizar.es). S. Lucia and B. Karg are with the Laboratory of Internet of Things for Smart Buildings, Technische Universität Berlin, Germany (e-mail: sergio.lucia@tu-berlin.de).

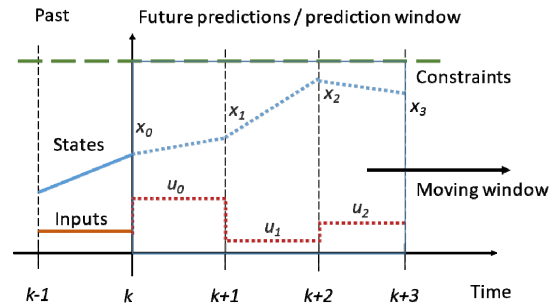

Fig. 1. Schematic description of Model Predictive Control. A mathematical model is used to predict the future states of a dynamic system and compute a sequence of optimal control inputs by solving an optimization problem.

optimal control inputs by solving online a numerical optimization problem that minimizes a chosen control goal (see Fig. 1).

MPC has been applied to power electronics converters [5, 17], specially to control electrical machines and drives. However, among the many applications presented in literature [18], little research has been published concerning resonant power converters.

Because of the computational complexity of MPC, most research and applications of MPC for power electronics has focused on strongly simplified formulations. The most commonly used [17], is the Finite Control Set (FCS) approach, which chooses a control action between a discrete amount of possible switching strategies by testing all possibilities. However, considering long horizons or variable switching frequencies is more challenging. Other approaches such as Explicit MPC [17] implement MPC as a look up table and are suitable for small linear systems. While there exist improvements for both the FCS and advanced Explicit MPC methods, it is still very challenging to design a controller that can consider simultaneously nonlinearities, long horizons, varying switching frequencies, time-varying constraints and can be easily implemented. The motivation of this work is to propose an approach that can achieve all these goals.

Instead of starting from strong simplifications, in this paper we design a rigorous nonlinear MPC problem. Motivated by the latest advances in the field of deep learning [19], [20], we use a deep neural network to directly learn the MPC solution. The online implementation of the controller is thus reduced to the evaluation of a neural network, which requires less resources than explicit MPC and can cope with nonlinear systems and long horizons, as recently shown in [21, 22]. The use of machine learning to approximate MPC controllers has also been recently proposed in [23], [24], [25], [26]. With respect to [23] and [25], our approach considers nonlinear systems. In contrast to [24] and [26], we use deep neural networks to achieve an improved accuracy and we include a detailed discussion on data generation. Further novel contributions that were not considered in previous results include an FPGA implementation, as well as Hardware-in-the-loop (HIL) results and the consideration of the tracking performance of the proposed controller under uncertainty of model parameters.

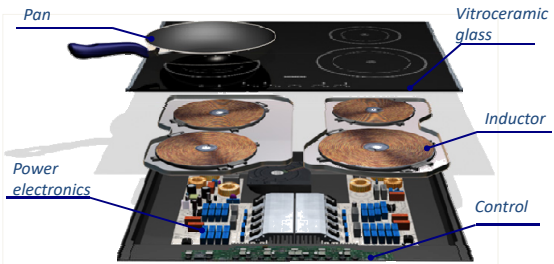

Fig. 2. Induction heating home appliance.

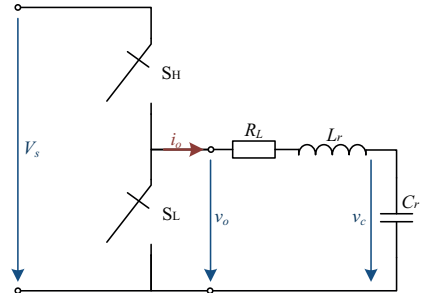

Fig. 3. Half-bridge series resonant converter.

The aim of this paper is to propose a new approach to design advanced MPC controllers to control resonant power converters. More specifically, the proposed MPC scheme will be applied to the control of a resonant inverter for induction heating (IH) applications. This application has a wide range of operating conditions which severely affect the resonant converter tuning, opening a wide field for control optimization using MPC. In the past, few research studies have dealt with IH applications. As far as the authors know, only in [27] a model predictive controller is proposed, with limited insight into the controller itself. Further contributions of this paper include a novel time transformation that enables a simple MPC formulation to explicitly handle important constraints such as zero-voltage switching. This paper extends the results presented in [28] by including novel contributions as the approximation of the MPC control law based on deep learning and incorporating FPGA implementation and HIL simulations. This implementation is intended for a final cost-effective ASIC implementation for domestic induction heating systems.

The remainder of this paper is organized as follows. Section II reviews IH technology and Section III covers model predictive control and its application to the half-bridge series resonant converter. Section IV and V summarizes the main neural network and FPGA implementation results, respectively. Finally, Section VI summarizes the main conclusions of this paper.

## II. RESONANT POWER CONVERSION FOR IH SYSTEMS

Induction heating [29] has become the leading heating technology due to its benefits in terms of efficiency, safety and performance. The applications of IH have been extended in recent years to cover a wide range of industrial, domestic and medical applications. Among these, domestic IH appliances are a relevant example where cookers (Fig. 2) ranging from simple burners to complex multiload systems [30, 31] are being developed.

Resonant power conversion plays a key role in such systems to create an alternating voltage to supply the induction coil, typically in the 20 kHz to 100 kHz frequency range. Usually, domestic IH systems implement the series-resonant half-bridge converter (Fig. 3) due to its good balance between cost and performance. The half-bridge series resonant inverter is composed of two switching devices, $S_H$ and $S_L$, typically implemented using IGTBs, and the resonant tank. The resonant tank is composed of the equivalent impedance of the coil-pot system, $R_L$, $L_r$, and the resonant capacitor $C_r$. The former parameters depend on the IH load selected by the user and its magnetic coupling, whereas the latter can be considered constant. Thus, the differential equation system that describes the dynamics of the resonant tank is defined by the current in the inductor, $i_o$,

$$R_L i_o(t) + L_r \frac{di_o(t)}{dt} + \frac{1}{C_r}\int i_o(t)dt = v_o(t). \qquad (1)$$

The modulation strategies typically used to control the output power include square wave modulation and asymmetrical duty cycle [32]. The former is based on modifying the switching frequency and the latter is based on modifying both, duty cycle $D$ and switching frequency $f_{sw}$. The maximum output power is obtained at the resonant frequency, defined as $f_o = 1/2\pi\sqrt{L_r C_r}$. In order to achieve zero voltage switching to reduce switching losses, $f_{sw}$ and $D$ are increased to reduce the output power. If the required output power is too low, typically pulse density modulation (PDM) is applied in order to limit the maximum $f_{sw}$.

The control strategy of the series resonant inverter for IH applications is usually constrained by several aspects inherent to the application, being the most important the following:

*Load variability*: The IH load, geometry and coupling, and required output power vary within a wide range, changing the resonant tank and, consequently, making the design of the controller a challenging task [33]. Classical PI controllers must be tuned to respond to fixed plant characteristics, whereas MPC controllers can be constrained to actuate in the desired manner under a wider range of operating conditions.

*Soft-switching*: In order to optimize the converter efficiency and safe operation, it is essential to guarantee soft-switching in the whole operational range [32], i.e. achieve ZVS conditions. This is specially challenging when using linear controllers, where some external constraints are difficult to be taken into account and may lead to unexpected performance and/or instability. MPC controllers allow to introduce these conditions, based on specific circuit states at the switching moment, as an additional constraint, whereas it is not feasible to use this approach with classical linear controllers.

*Multi-load control*: Induction heating systems are usually composed by several loads that may operate simultaneously. Under these circumstances, it is essential to guarantee that they operate synchronously, and no acoustic noise is generated due to switching frequencies intermodulation. In this context, MPC offers powerful tools for control optimization where a wide range of constrains, state space variables and control actions are present, optimizing the overall switching performance.

*Special operating conditions*: Under some conditions, the controller must be adapted due to excessive temperature, low electromagnetic coupling, maximum power consumed by several loads or limiting flicker emissions. All these specific issues are usually addressed by saturating or introducing exceptions in classical controllers. However, this can lead to unpredictable operation or put at risk the converter safety and performance. MPC provides a powerful control environment where all these restrictions can be taken into account without degrading the controller performance and reliability.

All these constrains severely limit the control possibilities of the converters and possess a significant challenge to design the control strategy. Moreover, as it has been discussed, all these constraints are difficult to be implemented using classical linear control schemes, increasing the complexity of the controller and limiting its scalability. For these reasons, MPC is explored as an effective control scheme to provide improved control and versatility for resonant power converters applied to IH.

### III. Model Predictive control

#### A. General formulation

Model predictive control uses the mathematical model of a dynamic system to predict its future behavior and compute a sequence of optimal control inputs that satisfies given constraints and optimizes a desired performance measure.

The model of a dynamic system can be written in a discrete-time setting as:

$$x_{k+1} = f(x_k, u_k), \qquad (2)$$

where $x_k \in \mathbb{R}^n$ denotes the state of the system at time $k$ and $u_k \in \mathbb{R}^m$ the control inputs. One of the main strengths of MPC is that it can consider nonlinear systems, general constraints, as well as control tasks different than classical setpoint tracking.

TABLE I
SUMMARY OF THE MODEL STATE AND INPUT VARIABLES

| Variables | Description |
|---|---|
| $i_o$ | State variable output current: Partly defines the second-order resonant circuit. It is proportional to the created magnetic field and, therefore, the induced heating. |
| $v_c$ | State variable resonant capacitor voltage: Completes the definition of the second-order resonant circuit. |
| $f_{sw}$ | Input variable switching frequency of the inverter devices: In a series resonant circuit, it is used to control the output power and is, typically, equal or higher than the resonant frequency to achieve soft-switching conditions. |
| $D$ | Input variable duty ratio of the half-bridge inverter devices: is the activation ratio of the switching devices during a switching period. In the proposed half-bridge configuration, the output voltage/power is maximum at $D = 0.5$, and they decrease when $D$ increases/decreases with symmetric behavior. |

The central challenge, especially for applications in power electronics, is that the implementation of MPC requires the online solution of the following optimization problem at each sampling time:

$$\begin{aligned} \underset{u_k}{\text{minimize}} & \quad \sum_{k=0}^{N-1} J(x_k, u_k) \\ \text{subject to} & \quad x_{k+1} = f(x_k, u_k), x_0 = \hat{x} \\ & \quad g(x_k, u_k) \leq 0 \end{aligned} \qquad (3)$$

where, $g$ describes general constraints and $\hat{x}$ is the current state of the system, which needs to be measured or estimated every time a new control should be computed. The performance metric, usually called cost function, is denoted by $J$ and $N$ is the prediction horizon.

In model predictive control, from the computed sequence of optimal control inputs ($u_0$, …, $u_{N-1}$) only the first element ($u_0$) is applied to the system and the same optimization problem is solved at the next step with the new state of the system. For applications that require sampling times in the range of microseconds, solving (2) in real-time can be very challenging despite of the recent progress in hardware, algorithms and tailored implementations [34] that have enabled the application of MPC to fast systems.

In the field of power electronics, most implementations of MPC [17] are based on a discretization of all control possibilities (FCS-MPC). A different possibility is the use of explicit MPC [35] which precomputes the solution of the optimization problem for all possible states, so that the online implementation reduces to a search algorithm.

Explicit MPC has been very successful for small linear systems but its implementation is more complex in the nonlinear case and for long horizons, because the amount of memory needed to store the controller grows exponentially with the prediction horizon and the number of constraints. A comparison of both continuous control set (CSS) and finite control set (FCS) MPC is presented in [36].

We propose in this work the use of deep learning to *learn*

the optimal control policy, as done in approximate explicit MPC [37], by means of a deep neural network (DNN) that can be easily implemented on an FPGA. This approach opens the door for the implementation of complex nonlinear MPC schemes with constraints and long horizons in fast systems, which can be implemented in FPGAs with a small need of digital resources that could be massively deployed in most power converters. We do not consider theoretical guarantees of stability or recursive feasibility as done e.g. in [38].

### B. MPC for Resonant Inverters

The resonant tank described in (1) can be modeled using the state-space paradigm. Table I summarizes the main state and input variables used in this model. Considering two states $x = [i_o, v_c]$, the system dynamics can be written as:

$$\begin{aligned} \frac{di_o}{dt} &= \frac{1}{L_r}(v_o - R_L i_o - v_c), \\ \frac{dv_c}{dt} &= \frac{1}{C_r} i_o. \end{aligned} \quad (4)$$

The available control inputs are the switching frequency and the duty cycle $u = [f_{sw}, D]$, which determine the signal $v_o$. One of the main challenges for the control of the half-bridge series resonant tank is the switched nature of $v_o$ which makes the dynamic system (3) a hybrid system. The signal $v_o$ can be written as:

$$v_o = \begin{cases} V_s & \text{if } t \in (t_i, t_i + D/f_{sw}] \\ 0 & \text{if } t \in (t_i + D/f_{sw}, t_{i+1}] \end{cases} \quad (5)$$

Solving general MPC problems of switched systems is a very difficult problem that leads to difficult optimization problems, especially if non-linearities and time-dependent constraints are considered, as it occurs for the computation of the average power or the consideration of ZVS constraints.

Relying on the finite control set (FCS) MPC technique is difficult in this case because we need to include time-varying constraints and make use of varying frequencies and duty cycles to improve efficiency. Instead, we propose to consider the full nonlinear model, with long horizons if necessary, and with flexible frequency and duty cycle in an efficient formulation that facilitates the offline solution of a large number of MPC problems. To enable the real-time implementation of the scheme, we propose to generate large amounts of data by solving rigorous nonlinear MPC problems and learn the control law using deep learning as will be explained in Section IV.

The proposed nonlinear model predictive controller for induction heating has three key elements: a time transformation, the average power computation, and the consideration of ZVS constraints. We explain these components in the remainder of the section

### C. Time transformation

The switching frequency and the duty cycle determine the switching instants $t_k$ of the system defined in (3). The direct consideration of the switching instants can be modeled using integer variables. However, if combined with long prediction horizons and time-varying constraints, it results in very complex mixed-integer optimization problems. We propose to perform a double time transformation that transforms the mixed-integer problem into a standard nonlinear but continous optimization problem.

We propose a time-varying transformation that maps the original time $t$ to the new scaled time units $\tau$ as follows:

$$\tau = \begin{cases} \dfrac{f_{sw,i}}{D_i} t & \text{if } t \in (t_i, t_i + D_i/f_{sw,i}] \\ \dfrac{f_{sw,i}}{1-D_i} t & \text{if } t \in (t_i + D_i/f_{sw,i}, t_{i+1}] \end{cases} \quad (6)$$

The main advantage of this transformation is that the switching in the new time units occurs exactly at $\tau = 1, 2, 3, \ldots$ We make use of this property to discretize the differential equations (3) in the new time units with a discretization time of $\tau_s = 1$ units. The model that is used in the MPC formulation can be written as:

$$x_{k+1}^{[mpc]} = \begin{cases} \dfrac{D_k}{f_{sw,k}} f^{[on]}(x_k, u_k) & \text{if } k = 0, 2, 4, \ldots, N \\ \dfrac{D_k - 1}{f_{sw,k}} f^{[off]}(x_k, u_k) & \text{if } k = 1, 3, 5, \ldots, N-1. \end{cases} \quad (7)$$

Where $f^{[on]}$ and $f^{[off]}$ denote discretized version of the differential equations in (3) when $v_o = V_s$ and $v_o = 0$ respectively.

This reformulation greatly simplifies the MPC formulation, as it is possible to consider long prediction horizons $N$ without the use of integer variables. It is also possible to use varying switching frequencies and duty cycles for each switching interval.

### D. Average Power Computation

To discretize the time-transformed version of the dynamics (7) to be used in the prediction of the MPC, we use orthogonal collocation on finite elements [39]. Each control interval is divided into finite elements, on which the state trajectory is parametrized using Lagrange polynomials. This discretization scheme is often used in nonlinear model predictive control as it provides a superior accuracy compared to simpler Euler discretization schemes. Using orthogonal collocation, the states of the system can be computed at any point in time and are directly available at the collocation points, which are new optimization variables.

The average power, calculated for each switching interval, can be computed as:

$$P_o = f_{sw} \int_0^{1/f_{sw}} i_o(t) v_o(t) dt \quad (8)$$

To achieve a tractable online computation of the average power, we use the value of the states at the collocation points within each control interval. Because of the time transformation previously described, each switching interval is composed of two control intervals (one step $k$ is used for the ON semi-cycle and one for the OFF semi-cycle). The average power at each control interval can be thus calculated as:

$$P_{o,k} = \begin{cases} f_{\text{sw},k} \sum_{i=1}^{n_{\text{col}}} i_{o,k}^{[i]} v_{o,k}^{[i]} (\tau_k^{[i]} - \tau_k^{[i-1]}), & \text{if } k = 0, 2, \ldots, N \\ P_{o,k-1}, & \text{if } k = 1, 3, \ldots, N-1, \end{cases} \quad (9)$$

where $n_{\text{col}}$ denotes the number of collocation points in each control interval and the $i_{o,k}^{[i]}$ denotes the value of the current at control step $k$ for collocation point $i$. The integral is approximated by a simple quadrature in which $\tau_k^{[i]}$ denotes the time of collocation point $i$ in the control interval $k$.

This implies that the average power can be predicted in advance, even during transient behavior, regardless of the inputs used (frequency, duty cycle, etc …) provided that the model is correct. Several techniques can be used to obtain good estimations of $R_L$ and $L_r$, see [33, 40] for an overview of possible methods. We also show in Section V that it is possible to adapt the MPC formulation to deal with uncertainty in $R_L$ and $L_r$.

*E. ZVS Constraints*

Another strength of the proposed time transformation is that enforcing zero-voltage-switching reduces to imposing static constraints at the end of each control interval. In the series resonant converter under study, ZVS depends on the operating conditions and can be achieved by ensuring the correct output current sign during the turn-on transition [3, 41]. These constraints can be easily incorporated in any MPC framework. It is straightforward to introduce any additional constraints that might be required by a specific application, as the ones described in Section II.

The optimization problem that should be solved after each switching cycle, based on the current measurements of the states $i_o$, $v_o$, and the power setpoint $P_o^{\text{des}}$ is:

$$\begin{aligned}
\underset{f_{\text{sw},k}, D_k}{\text{minimize}} \quad & \sum_{k=0}^{N} (P_{o,k} - P_o^{\text{des}})^2 + \alpha f_{\text{sw},k} & (a) \\
\text{subject to} \quad & \text{model in (7)} & (b) \\
& 30 \leq f_{\text{sw},k} \leq 100 \text{ kHz}, & (c) \\
& 0.2 \leq D_k \leq 0.8, & (d) \quad (10)\\
& i_{o,k} \geq 0, \text{if } k = 1, 3, 5, \ldots, N-1, & (e) \\
& i_{o,k} \leq 0, \text{if } k = 0, 2, 4, \ldots, N, & (f) \\
& D_k = D_{k-1}, \text{if } k = 1, 3, 5, \ldots, N-1, & (g) \\
& f_{\text{sw},k} = f_{\text{sw},k-1}, \text{if } k = 1, 3, 5, \ldots, N-1. & (h)
\end{aligned}$$

The tuning parameter α in the cost function (10a) weights the relative importance between using a small frequency and tracking the desired power. It is usually desired to choose the minimum possible frequency, as this maximizes the efficiency [32]. Note that the model used in the MPC problem is a discretized version of the model described in (4) using orthogonal collocation on finite elements and including the double time transformation described in (6). The discretized equations and power computation are included in constraints (10b), and (10c, 10d) denote the input constraints. The ZVS constraints are enforced by (10e, 10f). The fact that the inputs cannot be changed within a switching interval is enforced by (10g, 10h).

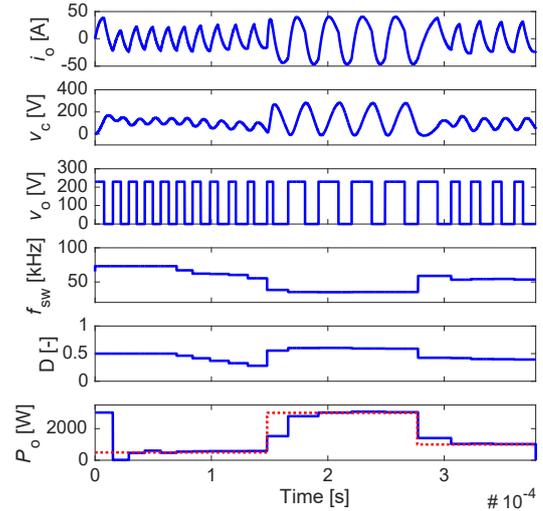

Fig. 4. Simulation results for the state variables, switching signal, control inputs and power. MPC controller when both $f_{\text{sw}}$ and $D$ are available as control inputs. The red dotted line denotes the desired power signal.

The optimization problem is nonconvex because of the power computation and because of the control-dependent time transformation.

The results presented here use the primal-dual interior point method described in [42] and implemented in IPOPT to solve all the optimization problems. Within the interior point algorithm, an exact Hessian is used. All the derivative information is computed using automatic differentiation via CasADi [43] and the MPC loop is implemented using the toolbox do-mpc [44]. The real system is simulated with the numerical integrator Sundials [45]. We use a prediction horizon of $N = 10$ steps, which equals 5 switching intervals. The length of the horizon is reasonable as it is the typical number of switching cycles that the system needs to achieve a new steady state. The model is discretized with orthogonal collocation based on Lagrange polynomials of degree 2, with 100 finite elements in each control interval to achieve an accurate average power computation. The tuning parameter that weights the importance in reducing the switching frequency is chosen as α =5e-8 in order to balance the relative contributions in the cost function.

The simulation results obtained by solving Problem (10) at each switching interval are shown in Fig. 4. The NMPC controller is switched on after 5 cycles and three different steps in the desired power are performed after 5 cycles each (500 W, 3000 W and 1000 W). The model parameters are $V_s$=230 V, $L_r$=19·10$^{-6}$ H, $R_L$=2.9 Ω, and $C_r$=1440 nF. It can be seen that the different power setpoints (red dashed line in the bottom plot of Fig. 4) can be tracked accurately while using both control inputs and respecting the ZSV constraints continuously, also during the transient behavior.

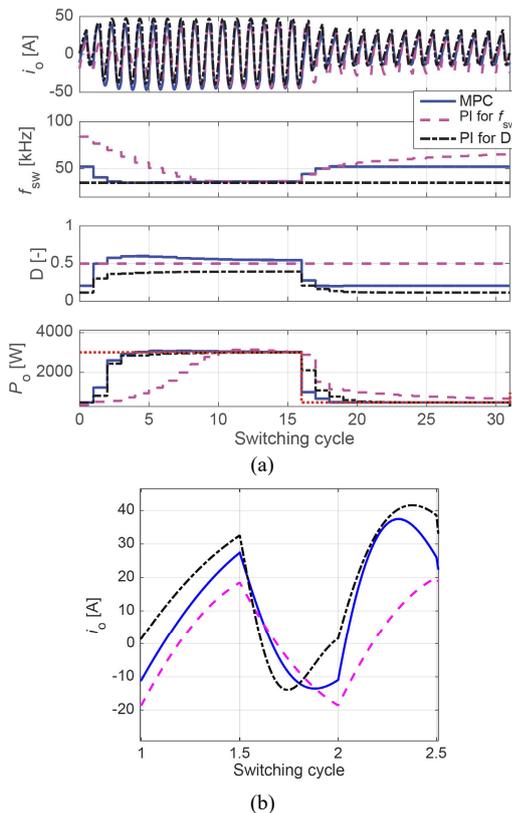

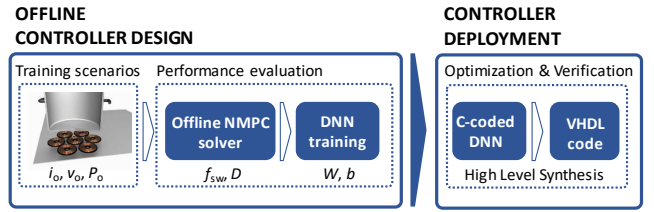

Fig. 6. Schematic representation of the controller design and deployment of the proposed deep learning-based NMPC solution.

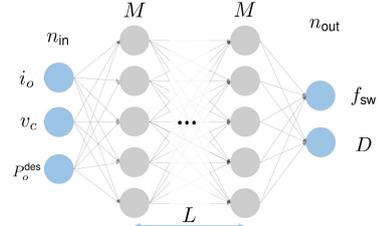

Fig. 7. Representation of a deep neural network with $L$ hidden layers with $M$ neurons each. The input layer and the output layer have $n_{in}$ and $n_{out}$ layers respectively.

Fig. 5. Comparison performance of the proposed MPC controller with PI implementations: (a) dynamic performance and (b) waveforms detail where $D$ PI do not achieve ZVS conditions.

To illustrate the advantages of the MPC controller, Fig. 5 (a) shows a comparison performance with classical PI controllers. It is important to note that solving problems for variable-plant MIMO systems would require otherwise using gain-schedulling techniques [46]. In this comparison, we consider a PI controller that uses the error in the power tracking to adapt the switching frequency $f_{sw}$ and another PI controller that uses the tracking error of the power to adapt the duty cycle $D$. The tuning of the PI controllers is performed via simulations and is chosen so that the fastest possible performance without overshoot in the delivered power is obtained. To facilitate the comparison between approaches, the x-axis represents switching cycles instead of time. It can be seen that the MPC approach is significantly faster than the PI controller for the switching frequency $f_{sw}$. While the PI controller for the duty cycle $D$ can obtain a similar performance when compared to MPC, it cannot incorporate constraints, resulting in violations of ZVS as it can be seen in Fig. 5 (b). Note that the MPC controller also enables the incorporation of any other constraint such that avoidance of acoustic noise due to frequency coupling between inductors. In same manner, several inductors could be controller simultaneously with the proposed MPC while solving problems for variable-plant MIMO systems would require otherwise using gain-scheduling techniques [45] that still cannot incorporate constraints.

The solution of each NMPC problems takes in average around 500 ms in a laptop computer (Macbook Pro, 3.5 GHz Intel Core i7 with 16 GB RAM). To obtain an implementable controller, we propose to approximate the solution of a rigorous NMPC formulation using deep neural networks instead of oversimplifying the problem formulation.

## IV. IMPLEMENTATION USING DEEP NEURAL NETWORKS

The embedded implementation of described controller is usually a challenging task due to real-time requirements and hardware restrictions. In this paper, an effective deep neural network implementation taking advantage of fast and parallel computing in FPGA will be proposed [6] (Fig. 6).

The approximation of the solution of MPC control laws has been usually studied under the headline of approximate explicit MPC [37]. While neural networks had been previously proposed to learn the mapping between states and control inputs [46], using deep neural networks for this purpose has been studied only very recently (see e.g. [21, 22]).

Using deep neural networks (with several hidden layers), instead of classical shallow ones (with only one hidden layer) has been one of the key developments that have enabled the practical successes of deep learning in previous years [19]. The latest theoretical advances [20] have shown that deep neural networks have greater expressivity than shallow networks, leading to better approximation capabilities of complex functions. This idea has been recently exploited in [22] to propose deep neural networks to approximate the function defined by the solution of the MPC problem. In that case, the inputs for the neural network are the current state of the system (and potentially also other parameters), and the outputs of the network are the optimal control inputs. Since the control strategy is defined via a neural network, our proposed approach could take advantage of future computing

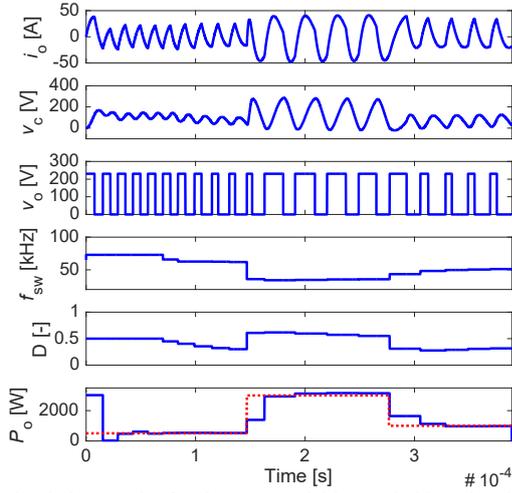

Fig. 8. Simulation results for the state variables, switching signal, control inputs and power when both $f_{sw}$ and $D$ are available control inputs computed via deep learning approximation of the MPC control law. The red dotted line denotes the desired power signal.

TABLE II
COMPARISON OF THE AVERAGE TRACKING PERFORMANCE INCLUDING TRANSIENT BEHAVIOR AND ZVS CONSTRAINT VIOLATIONS FOR 100 DIFFERENT RANDOM SETPOINT CHANGES WITH R AND L EXACTLY KNOWN

| Controller | Average tracking error [W/cycle] | Average ZVS violation [% of cycles] |
|---|---|---|
| Exact NMPC | 5.6939 | 0 |
| *Deep learning-based NMPC* | 6.2551 | 0 |

hardware tailored for the fast and efficient computation of machine learning tasks.

A deep neural network (Fig. 7) is a sequence of layers of neurons that determines a function $\mathcal{N}$ with inputs $z \in \mathbb{R}^{n_{in}}$ that can be defined as:

$$\mathcal{N}(z) = f_{L+1} \circ g_L \circ f_L \circ \ldots \circ g_1 \circ f_1(z), \quad (11)$$

where $L$ is the number of hidden layers and $M$ the number of neurons in each layer and $\circ$ denotes function composition. Each hidden layer consists of an affine function

$$f_l(\xi_{l-1}) = W_l \xi_{l-1} + b_l, \quad (12)$$

where $\xi_{l-1} \in \mathbb{R}^M$ is the output of the previous layer and $\xi_0 = z$. The matrices $W_l$ and the vectors $b_l$ are called the weights and biases at layer $l$ and their value is determined by training the neural network with known input-output pairs.
The second component of the neural network is the nonlinear function $g_l(\cdot)$, which is called activation function. Usual choices for $g_l(\cdot)$ include the rectifier linear unit (ReLU), the sigmoidal function or the *tanh* function.

Neural networks are *trained* offline, which means that the optimal values of the weights $W^{[l]}$ and biases $b^{[l]}$ are computed by minimizing the mean squared error:

$$\underset{W_l, b_l}{\text{minimize}} \quad \frac{1}{N_s} \sum_{s=1}^{N_s} (u^*(\hat{x}_s, P_o^{des}) - \mathcal{N}(\hat{x}_s, P_o^{des}))^2, \quad (13)$$

where $u^*$ denotes the optimal solution obtained when solving the MPC Problem (10) which is a function of the current state and the current setpoint. The inputs of the neural network $\mathcal{N}$ are also the current states and the current power setpoint and $N_s$ is the number of samples that are used for the training.

The training is done with Tensorflow [47] via Keras [48], using the optimizer ADAM [49], which is a modification of stochastic gradient descent and all computations are performed with 32-bit floating point number representation.

Fig. 8 shows the results of the deep-learning based controller for the same control task as the one shown in Fig. 4. It can be seen that the desired power setpoint is successfully tracked and the ZVS constraints are not violated. Both trajectories, for the exact solution of the NMPC and the proposed deep learning-based NMPC are very similar.

Table II presents a systematic evaluation of the performance of the proposed controller. We run 100 different simulations (similar to the ones shown in Fig. 4) for randomly chosen power setpoints. The different power setpoints are sampled between 500 W and 3000 W following a uniform distribution. As shown in Table II, the proposed controller achieves almost the same tracking performance than the exact solution of the NMPC problem and does not result in any violation of the ZVS constraints.

These promising results suggest that using the proposed method, it is possible to approximate very accurately a high-performance NMPC solution using a deep neural network that can be easily deployed on embedded hardware, as we show in the next sections.

## V. DATA GENERATION AND ROBUSTNESS ANALYSIS

The number of data pairs that are used to train the neural network and the way that such data pairs are generated has a significant influence on the approximation quality that can be achieved.

We chose a neural network with $L=5$ hidden layers and $M=10$ neurons per layer for all cases. We did not observe significant improvements for larger networks. The training data is composed of several input-output pairs. The input includes the power setpoint $P_o^{des}$, as well as the current measurements of $i_o$ and $v_c$ and the output includes the optimal inputs $f_{sw}$ and $D$ as illustrated in Fig. 7. For a given input, the output can be automatically generated by solving the MPC problem (10).

The input for which the corresponding outputs are calculated can be sampled randomly from the allowed state space in which the system can be operated. We sample this space using a uniform distribution for each one of the three inputs with the ranges $P_o^{des} \in [0, 3000]\,\text{W}$, $i_o \in [-150, 150]\,\text{A}$ and $v_c \in [-2000, 2000]\,\text{V}$.

If the optimization problem for the chosen initial state is unfeasible, the data point is discarded. The training error for different amounts of data points can be seen in Fig. 9.

An alternative approach to generate training data is to use an initial condition from the same allowed space but then simulate the closed-loop behavior using the MPC controller for a given amount of simulation steps. In this way, we generate training samples around optimal trajectories, which is where the system is usually operated and we avoid the use of states which might never be reached by the system. Using this kind of data generation, the training error can be reduced by an order of magnitude (see Fig. 9) and that is why this method is used to generate the data for the considered example. As expected, increasing the amount of data also leads to an improvement of the training error in any case. The validation error on previously unknown trajectories has a similar behaviour and it is omitted for brevity. Even if the network is trained on optimal trajectories, the obtained controller can be robust to different model parameters as we illustrate in the following.

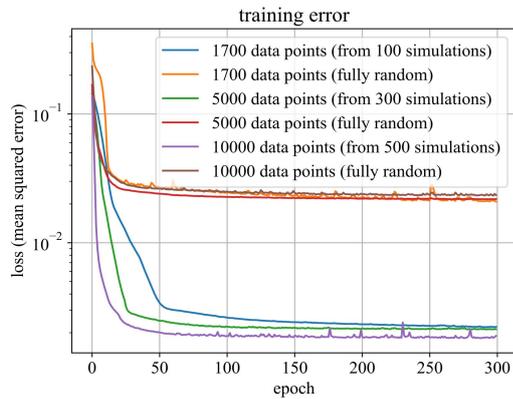

Fig. 9. Training error along the training iterations (epochs) for the same neural network structure when different amount of training data is used and different methods for the generation of the training data.

Ensuring the robustness of the controller with respect to model parameters is important since they are usually not exactly known or can vary rapidly. Achieving zero steady state error in NMPC under the presence of uncertainty can be often achieved by updating the desired power setpoint with the current power using a simple linear update rule. That is, at each sampling time the desired power used in the cost function of the MPC is updated following:

$$P_o^{\text{des}} = P_o^{\text{des}} + K(P_o^{\text{des, orig}} - P_o^{\text{des, meas}}), \qquad (14)$$

where $P_o^{\text{des, orig}}$ is the original desired power setpoint, $P_o^{\text{des, meas}}$ is the currently measured power, and $K = 0.8$ is a tuning parameter that can be chosen to balance the speed of the corrective actions. If the correction term (14) is not used, the formulation is the same as the one in (10), where $P_o^{\text{des}} = P_o^{\text{des,orig}}$ is constant and $P_o^{\text{des,meas}}$ is not used. An additional advantage of the proposed deep learning-based MPC controller is that exactly the same correction rule can be applied. In this case, the updated power setpoint is chosen as input to the neural network. As a result, steady state accuracy is achieved even under the presence of model errors. Fig. 10 shows the simulation results for the proposed approximated NMPC controller with a 15 % error in the model parameters $R_L$ and $L_r$. The correction strategy is activated at the middle of the simulation time and it can be seen that the correction leads to no tracking error despite of the uncertain parameters. We also performed 100 simulations for random values of $R_L$ and $L_r$ within a 15 % error interval with respect to the ones that were used for training. In this figure, ZVS operation can be also observed, verified by the correct current sign during the turn-on transition. A summary of the results can be seen in Table III. The tracking error is again very similar to the case when the optimization problem is solved at each sampling time. Although ZVS violations are slightly smaller for the deep learning-based controller, the difference is negligible and it in both cases only minor ZVS violations occur. If strict ZVS is desired for each cycle, robust NMPC techniques could be used.

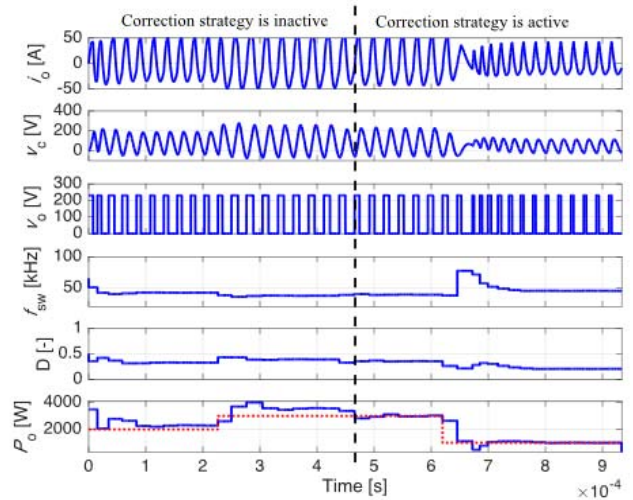

Fig. 10. Simulation results for the state variables, switching signal, control inputs and power when both $f_{sw}$ and $D$ are available control inputs computed via deep learning approximation of the MPC control law. The red dotted line denotes the desired power signal. There is a 15% error for the model parameters $R_L$ and $L_r$.

## VI. EMBEDDED IMPLEMENTATION RESULTS

All the embedded results presented in this section take advantage of the time transformation presented in (7) and which enables the fast solution of the MPC problems that generate the training data for the neural networks, as described in Section IV and Section V. The main motivation to approximate the solution of MPC via a neural network is that it enables a very simple embedded implementation of a complex controller, as explained in the remainder of the section.

### A. FPGA Implementation

The proposed MPC controller using deep neural networks is implemented using FPGA technology to take advantage of fast and parallel processing. Moreover, recent advances in FPGA technology [6] have led to cost-effective devices with an increasing amount of digital resources, enabling high-performance and competitive implementations.

Neural networks can be implemented with low precision arithmetic without degrading performance [50]. In order to quantize and normalize coefficients, in this work we use Ristretto [51], an automated Neural Networks approximation tool which condenses 32-bit floating point networks. Ristretto is an extension of Caffe [52] and allows to test, train and fine-tune networks with limited numerical precision.

In order to explore the design space and obtain an optimized implementation, automatic code generation using high level synthesis (HLS) [9, 53] has been applied. The implementation has been optimized to achieve the desired timing performance while minimizing the required digital HW requirements. As it is shown in Fig. 11, each neuron in each layer is calculated sequentially in order to minimize HW resources while achieving the required timing performance (see Table IV). By using the following HLS directive, the maximum number of multipliers is limited to 10:

HLS ALLOCATION instances=mul limit=10 operation

After that, the same digital HW is used to calculate the next layer, improving the final implementation. This is done by using directive:

HLS ALLOCATION instances=fcc limit=1 function

Finally, in order to optimize usage, each layer coefficients are stored in the same memory using dual-port memories. By doing that, only 5 memories are required to feed data to the 10 multipliers. This is implemented by using:

HLS array_map variable=weights_1 instance=weights horizontal
HLS array_map variable=weights_2 instance=weights horizontal

Table IV shows the performance for floating point and 16 bits implementation. The fixed-point version reaches the target latency of 1 µs while keeping the digital HW resources to low levels. The error with respect to the floating-point implementation is less than 5% in the worst case.

### B.A. HIL results

Simulation and verification of induction heating systems are challenging due to the high variety of operating conditions including output power levels and different resonant tanks determined by the pot materials, geometry and temperature. These parameters are usually randomly chosen by the user and, consequently, cannot be predicted. Furthermore, new materials and geometries usually appear in the market to be used with deployed appliances that should be operative for years. In this context, HIL has proven to be a very effective tool to simulate and test proposed controllers and control strategies under a wide variety of conditions [10, 33, 54] in a safer and faster way and, consequently, it will be used in this paper.

The proposed HIL platform used for domestic induction heating verification is depicted in Fig. 12. It consists on an FPGA-based architecture that runs both ad-hoc digital hardware and software running in the dual-core ARM processor available. Linux OS is used to generate automatically the required test-benches as well as to perform results verifications and communications with a user interface. A plant model of the series resonant inverter is also implemented, which operates in real time, taking as inputs the control variables and the operating parameters, i.e. resonant tank, and generating as outputs the state variables. The proposed deep learning-based MPC controller is implemented

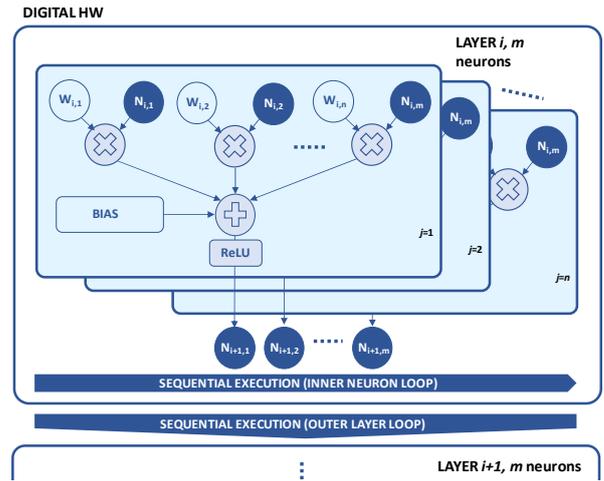

Fig. 11. Optimized arithmetic FPGA implementation.

TABLE III
COMPARISON OF THE AVERAGE TRACKING PERFORMANCE INCLUDING TRANSIENT BEHAVIOR AND ZVS CONSTRAINT VIOLATIONS FOR 100 DIFFERENT RANDOM SETPOINT CHANGES WITH 15 % ERROR IN R AND L

| Controller | Average tracking error [W/cycle] | Average ZVS violation [% of cycles] |
|---|---|---|
| Exact NMPC | 10.845 | 0.3922 |
| Deep learning-based NMPC | 12.243 | 0.0784 |

TABLE IV
FPGA IMPLEMENTATIONS SUMMARY

| Item | 16-bit implementation | Floating-point implementation |
|---|---|---|
| Latency | 1 µs | 2.1 µs |
| LUTs | 1382 | 12851 |
| FFs | 888 | 7396 |
| DSPs | 10 | 50 |

in the FPGA following the implementations detailed at the beginning of this Section.

Fig. 13 shows an example of a post-layout simulation with 16-bit arithmetic. In this figure, the current through the resonant tank, the target output power and the controller error are represented. As it can be seen, the controller works as expected, achieving the desired output power while keeping ZVS constraints, which can be verified checking the proper current sign during the switching times. These results prove the feasibility of the proposed controller, as well as its FPGA-implementation using deep neural networks and high-level synthesis.

To further validate the robustness of the proposed controller, we performed several HiL experiments using the deep learning-based NMPC with the correction term (14) which are shown in Table V. In particular, we perform nine different experiments where three different steps in the desired power output $P_o^{des,orig}$ are performed. The results reported in Table V include the steady state accuracy (after 10 switching cycles) of the controller with 16-bit arithmetic as well as the number of ZVS violations for different values of the parameters $R_L$ and $L_r$. In all cases, the proposed controller achieves a very small steady state error for different values of the parameters thanks to the correction term (14).

## VII. CONCLUSIONS

Resonant power conversion enables the implementation of power converters with superior performance, efficiency and power density. Among the multiple applications, induction heating systems are a relevant example which outperforms other technologies due to its superior efficiency and performance. However, the accurate control of the resonant converter under a wide variety of constraints of different nature still remains a challenge. This paper has proposed how high-performance NMPC solutions can be easily and accurately approximated using deep neural networks. The resulting neural networks can be used in a High-Level Synthesis framework to obtain FPGA designs that enable the real time advanced control of power converters.

The proposed scheme has been detailed, and the controller performance has been successfully validated in simulations and in a hardware-in-the-loop setup. The proposed controller enables to take into account complex control constraints for IH systems as well as to obtain improve dynamics in a variable-plant MIMO system that would require otherwise more complex and difficult to tune control techniques. Furthermore, the proposed control architecture enables the implementation of a complex MPC controller in a simple manner in state-of-the-art ASIC implementations with significant benefits in terms of control constraints and dynamic performance. As a conclusion, the proposed approach opens the door to the application of complex NMPC-schemes for future higher-performance higher-complexity IH systems. Future work will study the scalability properties of the proposed approach and further comparisons with other approximate explicit MPC methods as done in the linear case in [22].

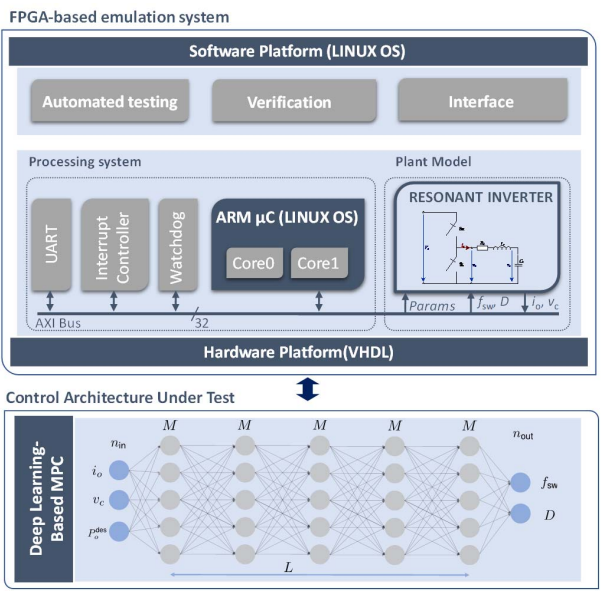

Fig. 12. Hardware-in-the-loop testbench.

TABLE V
STEADY STATE TRACKING ERRORS FOR THE HIL IMPLEMENTATION WITH DIFFERENT SETPOINTS AND DIFFERENT PARAMETER VALUES

| R error | L error | Steady state error (mW) for different setpoints | | | Average ZVS violation [%] |
|---|---|---|---|---|---|
| | | 1000 W | 2000 W | 3000 W | |
| 0% | 0% | 41.74 | 1.554 | 2.747 | 0 |
| 0% | -15% | 6.903 | 2.149 | 0.699 | 0 |
| 0% | +15% | 96.22 | 448.1 | 737.7 | 0 |
| -15% | 0% | 10.39 | 2.023 | 97.46 | 0 |
| -15% | -15% | 8.632 | 3.964 | 5.97 | 0 |
| -15% | +15% | 223.0 | 449.3 | 54.19 | 0 |
| +15% | 0% | 48.46 | 22.68 | 125.1 | 0 |
| +15% | -15% | 2.234 | 3.282 | 3.693 | 0 |
| +15% | +15% | 49.78 | 116.9 | 108.3 | 0 |

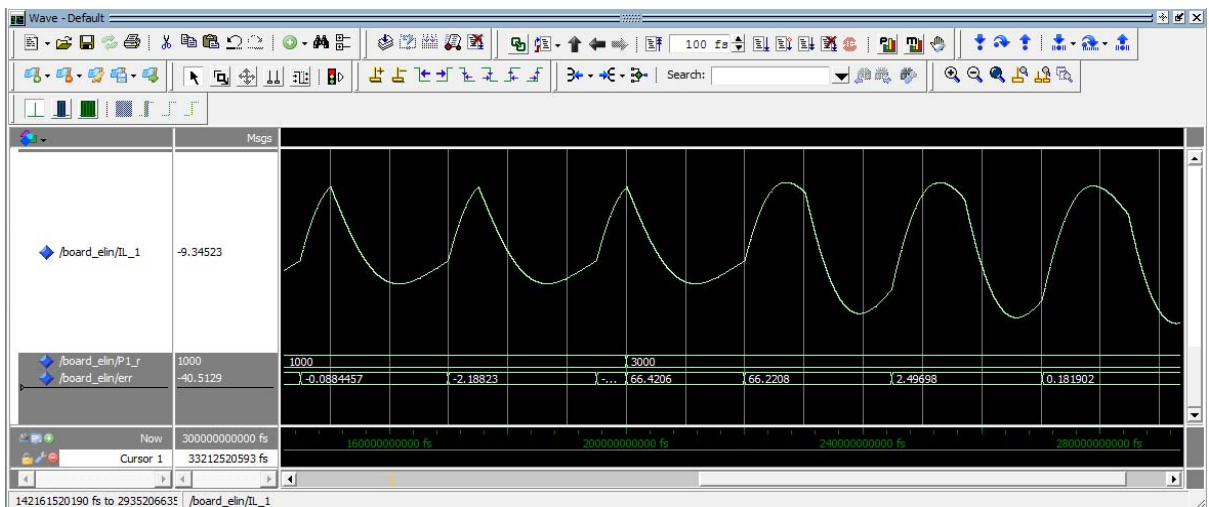

Fig. 13. Post-layout simulation using 16-bit arithmetic.

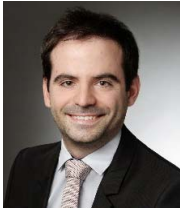
**Sergio Lucia** (M'16) received the M.Sc. degree in electrical engineering from the University of Zaragoza, Zaragoza, Spain, in 2010, and the Dr. Ing. degree in optimization and automatic control from the Technical University of Dortmund, Dortmund, Germany, in 2014. He joined the Otto-von-Guericke Universitat Magdeburg and visited the Massachusetts Institute of Technology as a Postdoctoral Fellow.
Since May 2017, he has been an Assistant Professor and Chair with the Laboratory of Internet of Things for Smart Buildings, Technische Universitat Berlin, Berlin, Germany, and with Einstein Center Digital Future, Berlin. His research interests include decision-making under uncertainty, distributed control, as well as the interplay between machine learning techniques and control theory. Dr. Lucia is currently Associate Editor of the Journal of Process Control.

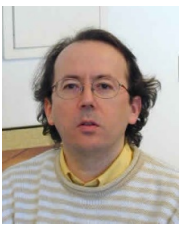
**Denis Navarro** received the M.Sc. degree in Microelectronics from the University of Montpellier, France, and the Ph.D. degree from the University of Zaragoza in 1987 and 1992, respectively.
Since September 1988, he has been with the Department of Electronic Engineering and Communications at the University of Zaragoza, where he is a Professor. His current research interests include CAD for VLSI, low power ASIC design, and modulation techniques for power converters. He is involved in the implementation of new applications of integrated circuits. In 1993 Dr. Navarro designed the first SPARC® microprocessor in Europe.
Dr. Navarro is a member of the Aragon Institute for Engineering Research (I3A).

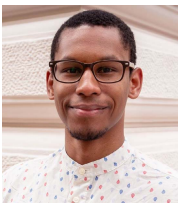
**Benjamin Karg** was born in Burglengenfeld, Germany, in 1992. He received the B.Eng. degree in mechanical engineering from Ostbayerische Technische Hochschule Regensburg, Regensburg, Germany, in 2015, and his M.Sc. degree in systems engineering and engineering cybernetics from Otto-von-Guericke Universität, Magdeburg, Saxony-Anhalt, Germany, in 2017. He currently works as a research assistant at the laboratory "Internet of Things for Smart Buildings", Technische Universität Berlin, Germany, to pursue his PhD. He is also member of the Einstein Center for Digital Future. His research is focused on control engineering, artificial intelligence and edge computing for IoT-enabled cyber-physical systems.

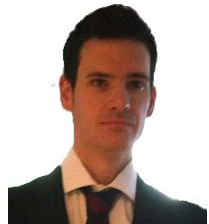
**Hector Sarnago** (S'09 M'15 SM'19) received the M.Sc. degree in Electrical Engineering and the Ph.D. degree in Power Electronics from the University of Zaragoza, Spain, in 2010 and 2013, respectively. Currently, he is a senior post-doc researcher in the the Department of Electronic Engineering and Communications at the University of Zaragoza, Spain. His main research interests include resonant converters and digital control for induction heating applications.
Dr. Sarnago is a member of the Aragon Institute for Engineering Research (I3A).

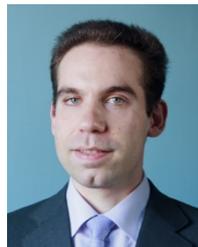
**Óscar Lucía** (S'04, M'11, SM'14) received the M.Sc. and Ph.D. degrees (with honors) in Electrical Engineering from the University of Zaragoza, Spain, in 2006 and 2010, respectively.
During 2006 and 2007 he held a research internship at the Bosch and Siemens Home Appliances Group. Since 2008, he has been with the Department of Electronic Engineering and Communications at the University of Zaragoza, Spain, where he is currently an Associate Professor. During part of 2009 and 2012, he was a visiting scholar at the Center of Power Electronics Systems (CPES), Virginia Tech. His main research interests include resonant power conversion, wide-bandgap devices, and digital control, mainly applied to contactless energy transfer, induction heating, electric vehicles, and biomedical applications. In these topics, he has published more than 70 international journal papers and 150 conference papers, and he has filed more than 40 international patents.
Dr. Lucía is a Senior Member of the IEEE and an active member of the Power Electronics (PELS) and Industrial Electronics (IES) societies. He was a Guest Associate Editor of the IEEE Transactions on Industrial Electronics and the IEEE Journal of Emerging and Selected Topics in Power Electronics in 2013 and 2015, respectively. Currently, he is an Associate Editor of the IEEE Transactions on Industrial Electronics, IEEE Open Journal of Industrial Electronics, and IEEE Transactions on Power Electronics. Dr. Lucía is a member of the Aragon Institute for Engineering Research (I3A).